\documentclass[12pt,a4]{amsart}
\usepackage{amssymb}
\usepackage{amsmath}
\usepackage{amscd}
\usepackage{graphicx}

\usepackage[utf8]{inputenc}
\usepackage[usenames,dvipsnames,pdftex]{xcolor}
\usepackage[upright]{fourier}
\usepackage{tkz-graph}
\usetikzlibrary{arrows}
\usetikzlibrary{decorations.markings}

\newtheorem{proposition}{Proposition}[section]
\newtheorem{theorem}[proposition]{Theorem}
\newtheorem{lem}[proposition]{Lemma}
\newtheorem{cor}[proposition]{Corollary}

\theoremstyle{remark} 
\theoremstyle{definition}
\newtheorem{exam}[proposition]{Example}
\newtheorem{rem}[proposition]{Remark}

\numberwithin{equation}{section}

\newcommand{\Z} {\mathbb{Z}}
\newcommand{\N} {\mathbb{N}}
\newcommand{\R} {\mathbb{R}}
\newcommand{\C} {\mathbb{C}}
\newcommand{\Q} {\mathbb{Q}}
\newcommand{\A} {\mathcal{A}}

\newcommand{\X}{\mathcal{X}}
\newcommand{\Y}{\mathcal{Y}}
\DeclareMathOperator{\supp}{supp}
\DeclareMathOperator{\einf}{essinf}
\renewcommand{\L}{\mathcal{L}}

\begin{document}

\title{Rotational beta expansion: Ergodicity and Soficness}
\author{Shigeki Akiyama}
\author{Jonathan Caalim}
\thanks{The authors are supported by the
Japanese Society for the Promotion of Science (JSPS), Grant in aid
21540012. The second author expresses his deepest gratitude to the 
Hitachi Scholarship Foundation.
}
\email{akiyama@math.tsukuba.ac.jp \hspace{10mm} nathan.caalim@gmail.com
}
\address{Institute of Mathematics \&
Center for Integrated Research in Fundamental Science and Technology, 
University of Tsukuba, 1-1-1 Tennodai,
Tsukuba, Ibaraki, Japan (zip:350-8571)}
\address{Institute of Mathematics, University of the Philippines Diliman,
1101 Quezon City, Philippines}

\maketitle	

\begin{abstract}
We study a family of piecewise expanding maps on the plane, generated 
by composition of a rotation and an expansive similitude of expansion constant $\beta$. 
We give two constants $B_1$ and $B_2$ depending 
only on the fundamental domain that
if $\beta>B_1$ then the expanding map
has a unique absolutely continuous invariant probability measure, and if
$\beta>B_2$ then it is equivalent to $2$-dimensional Lebesgue
measure.
Restricting to a rotation generated by $q$-th root of unity $\zeta$
with all parameters in $\Q(\zeta,\beta)$, 
the map gives rise to
a sofic system 
when $\cos(2\pi/q) \in \Q(\beta)$ and $\beta$ is a Pisot number. It is also
shown that the condition $\cos(2\pi/q) \in \Q(\beta)$ is necessary by giving 
a family of non-sofic systems for $q=5$. 
\end{abstract}

\begin{section}{Introduction}

Let $1<\beta\in \R$ and $\zeta\in \C\setminus \R$ with $|\zeta|=1$. 
Fix $\xi,\eta_1,\eta_2 \in \C$ with $\eta_1/\eta_2 \not \in \R$.
Then  $\X=\{ \xi+x \eta_1+ y \eta_2 \ |\ x\in [0,1), y\in [0,1) \}$ is a 
fundamental domain of the lattice $\L$ generated by $\eta_1$ and $\eta_2$ in $\C$, i.e.,
$$
\C = \bigcup_{d\in \L} (\X + d)
$$
is a disjoint partition of $\C$.
Define
a map $T:\X \rightarrow \X$ by $T(z)= \beta \zeta z - d$
where $d=d(z)$ is the unique element in $\L$ 
satisfying $\beta \zeta z\in \X+d$.
Given a point $z$ in $\X$, we obtain an expansion
\begin{eqnarray*}
z &=& \frac{d_1}{\beta \zeta}+ \frac{T(z)}{\beta \zeta}\\
  &=& \frac{d_1}{\beta \zeta}+ \frac{d_2}{(\beta \zeta)^2} + \frac{T^2(z)}{(\beta \zeta)^2} \\
  &=& \sum_{n=1}^{\infty} \frac {d_n}{(\beta \zeta)^n}
\end{eqnarray*}
with $d_n=d(T^{n-1}(z))$. We call $T$ the {\it rotational beta transformation} 
and $d_1d_2...$ the {\it expansion} of $z$ with respect to $T$.
We note that
the map $T$ generalizes the notions of beta expansion \cite{Renyi:57,Parry:60,Ito-Takahashi:74}
and negative beta expansion \cite{Ito-Sadahiro:09, Liao-Steiner:12, Kalle:14}
in a natural dynamical manner to the complex plane $\C$. 
More number theoretical generalizations had been 
studied by means of numeration system in complex bases, e.g., 
\cite{Katai-Kovacs:81,Gilbert:81,
Akiyama-Brunotte-Pethoe-Thuswaldner:05,
Komornik-Loreti:07}.
Since $T$ is a piecewise expanding 
map, by a general theory developed in \cite{Keller:79, Keller:85, Gora-Boyarsky:89, Saussol:00, Tsujii:00, Buzzi-Keller:01, Tsujii:01},
there exists an invariant probability measure $\mu$ which is absolutely continuous to
the two-dimensional Lebesgue measure\footnote{For example, we can see this fact 
by Lemma 2.1 of \cite{Saussol:00} for some iterate of $T$.}.
The number of ergodic components is known to be finite
\cite{Keller:79, Gora-Boyarsky:89, Saussol:00}. An explicit upper bound 
in terms of the constants in a Lasota-Yorke type inequality
was given by Saussol \cite{Saussol:00}. However this bound may be 
large\footnote{Saussol \cite{Saussol:00} did not aim at giving a good bound of it, but
was interested in showing the finiteness of the number of components. 
Indeed, when we apply Lasota-Yorke type inequality, 
these two objectives (finiteness proof and minimizing the upper bound)
are in confrontation.}.
By using the special shape of the map $T$, we can show
that the number is one if $\beta$ is sufficiently large.
Define the {\it width} $w(\X)$ 
of $\X$ as 
$$
w(\X):=\min \left\{|\eta_1|,|\eta_2| \right\}\sin (\theta(\X)),
$$ 
where $\theta(\X)\in (0,\pi)$ is the angle between $\eta_1$ and $\eta_2$. Then
$w(\X)$ is the minimum height of the 
parallelogram
formed by $\X$.
Let $r(P)$ be the {\it covering radius} of a point set $P\subseteq \C$, 
i.e., $r(P)$ is the infimum of the positive real numbers $R$ 
such that every point in $\C$ is within distance $R$ of at least one point in $P$.
Let us define
$$
B_n=\max\left\{ \nu_n(\theta(\X)), \frac{2r(\L)}{w(\X)} \right\}\quad (n=1,2)
$$
with
$$
\nu_1(\theta(\X)):=\begin{cases} 
      2 & \mbox{ if } \frac{1}{2} < \tan \left(\frac{\theta(\X)}2\right)< 2\\
      \frac{1+|\cos\theta(\X)|}{2(\sin\theta(\X)+|\cos\theta(\X)|-1)}
      & \mbox{ otherwise } 
   \end{cases}
$$
and 
$$
\nu_2(\theta(\X)):=1+ \frac{\sqrt{2}}{\sin \theta(\X) \sqrt{
1+|\cos \theta(\X)|}}=1+\frac 1{\sin \theta(\X) 
\max\left\{\sin \frac{\theta(\X)}2,\cos\frac{\theta(\X)}2\right\}}.
$$
Note that $B_1$ and $B_2$ do not depend on $\xi$ and are
determined only by $\eta_1$ and $\eta_2$.
\begin{theorem}\label{Ergodic}
If 
$\beta > B_1$ then $(\X,T)$ has a unique absolutely continuous invariant 
probability measure $\mu$.
Moreover, if $\beta>B_2$
then $\mu$ is equivalent to the 2-dimensional Lebesgue measure restricted to $\X$.
\end{theorem}

One can confirm the inequality $B_1\le B_2$ in Figure \ref{bound}.
\begin{figure}[htp]
\caption{Comparison of $\nu_1(\theta(\X))$ and $\nu_2(\theta(\X))$}
	\centering
		\includegraphics[scale=0.7]{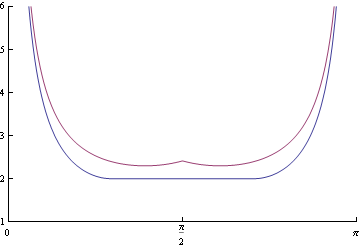}\label{bound}
\end{figure}
The uniqueness implies that $(\X,T)$ is ergodic with respect to $\mu$. 
In the last section, we give a rotational beta transformation where 
the number of ergodic components exceeds one, when $\beta$ is small 
(see Example \ref{NonErgode}).
It is an intriguing problem to improve the above bounds $B_1$ and $B_2$,
which may not be optimal, see Examples
\ref{Four}, \ref{Five} and \ref{Seven}.
Hereafter, ACIM stands for 
absolutely continuous invariant probability measure.

\begin{rem}
The covering radius $r(\L)$ is computed from the successive minima of $\L$, 
which are derived by the `homogeneous' continued fraction algorithm due to Gauss.
The term $2r(\L)/w(\X)$ in Theorem \ref{Ergodic} 
is expected to be replaced by a smaller one, 
since we may substitute $r(\L)$ with
$r(\L+T^{-n}(z))$ for a non negative integer $n$ and a point $z$ in $\X$
to obtain the same conclusion. See the proof in \S \ref{E}.
\end{rem}

\begin{rem}
The beta and negative beta transformations could be understood 
in a similar framework in 1-dimension by choosing 
$\zeta=\pm 1$ and $\X=[\xi, \xi+\eta)$
with $\L=\eta \Z$. In this case, $(\X,T)$ has a 
unique ACIM
with respect to the 1-dimensional Lebesgue measure. This result
follows from Li-Yorke \cite{Li-Yorke:78}
which reads that every support of an ACIM
contains at least one discontinuity in its interior, and the fact that
a neighborhood of each discontinuity of $T$ is
mapped similarly to neighborhoods of two end points of $\X$.
The problem of discontinuities becomes harder in dimension $>1$.
\end{rem}

Later on, we are interested in the 
associated symbolic dynamical system
over the alphabet $\mathcal{A} :=\{ d(z)\ |\ z\in \X\}$.
Let $\mathcal{A}^{\Z}$ (resp. $\mathcal{A}^*$) be the 
set of all bi-infinite (resp. finite) words over 
$\mathcal{A}$. 
We say $w\in \mathcal{A}^*$ is admissible
if $w$ appears in the expansion $d_1d_2...$ for some $z\in \X\setminus
\bigcup_{n=-\infty}^{\infty} T^n(\partial(\X))$\footnote{We exclude the null set $\bigcup_{n=-\infty}^{\infty} T^n(\partial(\X))$, i.e., the set of 
forward/backward discontinuities to concentrate on the essential part of the dynamics.}.
Let $$
\X_T:=\left.\left\{w=(w_j)\in \mathcal{A}^{\Z} \right| w_jw_{j+1}...w_k 
\mbox{ is admissible }\forall (j,k) \in \Z^2 \text{ with } j\le k \right\}
$$
which is compact by the product topology of $\A^{\Z}$. 
The symbolic dynamical system
associated to $T$ is the topological dynamics $(\X_T,s)$
given by the shift operator $s((w_j))=(w_{j+1})$.
We say $(\X_T,s)$ (or simply, $(\X,T)$) is sofic if there is a finite directed graph $G$ labeled
by $\mathcal{A}$ such that for each $w \in \X_T$, there exists a 
bi-infinite path in $G$ labeled $w$ and vice versa. 
Here is a characterization 
of sofic systems using the forward orbits of 
the discontinuities:

\begin{lem}\label{Sofic Definition}
The system $(\X,T)$ is sofic if and only if $\bigcup_{n=1}^{\infty} T^n(\partial(\X))$ is a finite union of segments.
\end{lem}

Here $\partial(\X)$ denotes the boundary of $\X$. Note that the two open segments in $\partial(\X)$, one from $\xi+\eta_1$ to
$\xi+\eta_1+\eta_2$ and the other from $\xi+\eta_2$ to $\xi+\eta_1+\eta_2$,
are outside of $\X$. 
For these segments, the images by $T$ are defined 
by an infinitesimal small perturbation, e.g., we take the image of the 
segment connecting $\xi+\eta_1(1-\varepsilon)$ and $\xi+\eta_1(1-\varepsilon)+\eta_2$ 
for a small positive $\varepsilon$. We prove this lemma in \S \ref{SD}.

From the above lemma, we see that for $(\X,T)$ to be sofic, 
the set of slopes of the discontinuous segments consisting
$\bigcup_{n=1}^{\infty} T^n(\partial(\X))$ 
must be finite. This means that
$\zeta$ must be a root of unity.
Hereafter, we assume that $\zeta$ is a $q$-th root of unity
with $q>2$
and $\xi, \eta_1, \eta_2\in \Q(\zeta,\beta)$ with $\eta_1/\eta_2\not \in \R$.
We let $\kappa(\xi+\eta_1 x+\eta_2 y)=\begin{pmatrix} x\\ y \end{pmatrix}$
be a bijection from $\X$ to $[0,1)^2$ and consider the analog of $T$ 
on $[0,1)^2$.

Since $\Q(\zeta,\beta)$ is quadratic over $\Q(\zeta+\zeta^{-1},\beta)$,
every element of $\Q(\zeta,\beta)$ is uniquely expressed as 
a linear combination of $\eta_1$ and $\eta_2$ over $\Q(\zeta+\zeta^{-1},\beta)$.
We find  $a_{jk}, b_j\in \Q(\zeta+\zeta^{-1},\beta)$ such that
$$
\zeta \begin{pmatrix} \eta_1\\ \eta_2 \end{pmatrix}
= \begin{pmatrix} a_{11}& a_{21} \\ a_{12}& a_{22}\end{pmatrix} 
\begin{pmatrix} \eta_1\\ \eta_2 \end{pmatrix}
$$
and
$$
(\beta \zeta-1) \xi = b_1\eta_1+b_2 \eta_2.
$$
Let $U$ be the map from $[0,1)^2$ to itself, which satisfies
$U \circ \kappa =\kappa \circ T$. We can write
$$
U\left(\begin{pmatrix}x\\y\end{pmatrix}\right)= 
\begin{pmatrix} \beta (a_{11} x+a_{12}y) + b_1 
-\lfloor \beta (a_{11} x+a_{12}y) + b_1 \rfloor\\ 
\beta (a_{21} x+a_{22}y) + b_2 
-\lfloor \beta (a_{21} x+a_{22}y) + b_2 \rfloor
\end{pmatrix}.
$$
This expression suggests an important role of the field $\Q(\zeta+\zeta^{-1})$
in our problem. In the following, we give a sufficient condition so that $(\X,T)$
is a sofic system.

\begin{theorem}
\label{Sofic}
Let $\zeta$ be a $q$-th root of unity $(q>2)$ and $\beta$ be a Pisot number. Let 
$\eta_1,\eta_2, \xi \in \Q(\zeta,\beta)$. 
If $\zeta+\zeta^{-1}\in \Q(\beta)$, then the system $(\X,T)$ is sofic. 
\end{theorem}

In proving this theorem, we give
an upper bound on the number of the {\it intercepts} 
of the segments in $\bigcup_{n=1}^{\infty} T^n(\partial(\X))$. 
The details will be given in \S \ref{S}. For $q=3,4,6$, since $2\cos(2\pi p/q)$ 
is an integer, we have the following result.

\begin{cor}
If $\zeta$ is a $3rd$, $4th$ or $6th$ root of unity, 
then the system $(\X,T)$ is sofic for any Pisot number $\beta$. 
\end{cor}

On the other hand, we can give a family of non-sofic systems when 
$\zeta+\zeta^{-1}\not \in \Q(\beta)$. From here on, $i$ denotes $\sqrt{-1}$.

\begin{theorem}
\label{Non-Sofic}
Let $\xi=0$, $\eta_1=1$ and $\eta_2=\zeta=\exp(2\pi i/5)$. 
If $\beta>2.90332$ such that
$\sqrt{5}\not \in \Q(\beta)$, then $(\X,T)$ is
not a sofic system.
\end{theorem}

Most of the large Pisot numbers satisfy the conditions of Theorem
\ref{Non-Sofic}, e.g., any integer greater than 2.
The proof of Theorem \ref{Non-Sofic} suggests that $(\X,T)$ rarely becomes
sofic for general $\beta$ and $\zeta$. 
Meanwhile, Example \ref{Four} shows that there are sofic rotational beta expansions 
beyond Theorem \ref{Sofic}. 
It is of interest to characterize such quintuples $(\beta, \zeta, \eta_1,\eta_2,\xi)$, 
giving an analogy of Parry numbers in $1$-dimensional beta expansion (cf.
\cite{Parry:60,Ito-Takahashi:74, Akiyama:15}).

\end{section}

\begin{section}{Proof of Theorem \ref{Ergodic}}\label{E}

Let $t$ be a positive real number. 
We denote by $B_{-t}(A)$ the set of points of $A$ which have 
distance at least $t$
from $\partial(A)$.
We shall study the $n$-th inverse image 
$T^{-n}(z)=\{z'\in \X| T^n(z')=z\}$ for $n\in \mathbb{N}$ and $z\in \X$. 
Put $r=r(\L), w=w(\X)$ and $\theta=\theta(\X)$. For $j=1,2$, set $\nu_j=\nu_j(\theta(\X))$.
First we claim that if $\beta > B_2$, 
then for all $z\in \X$, $\bigcup_{n=1}^{\infty} T^{-n}(z)$ is dense in $\X$.
Note that 
$T^{-1}(z)=\frac{1}{\beta \zeta} \left(\left(z+\mathcal{L}\right)\cap \beta\zeta \X\right)$. 

\begin{figure}[htp]
\caption{$B_{-\frac{r}{\beta}}\left( \X\right)+\L$}
	\centering
		\includegraphics[scale=0.10]{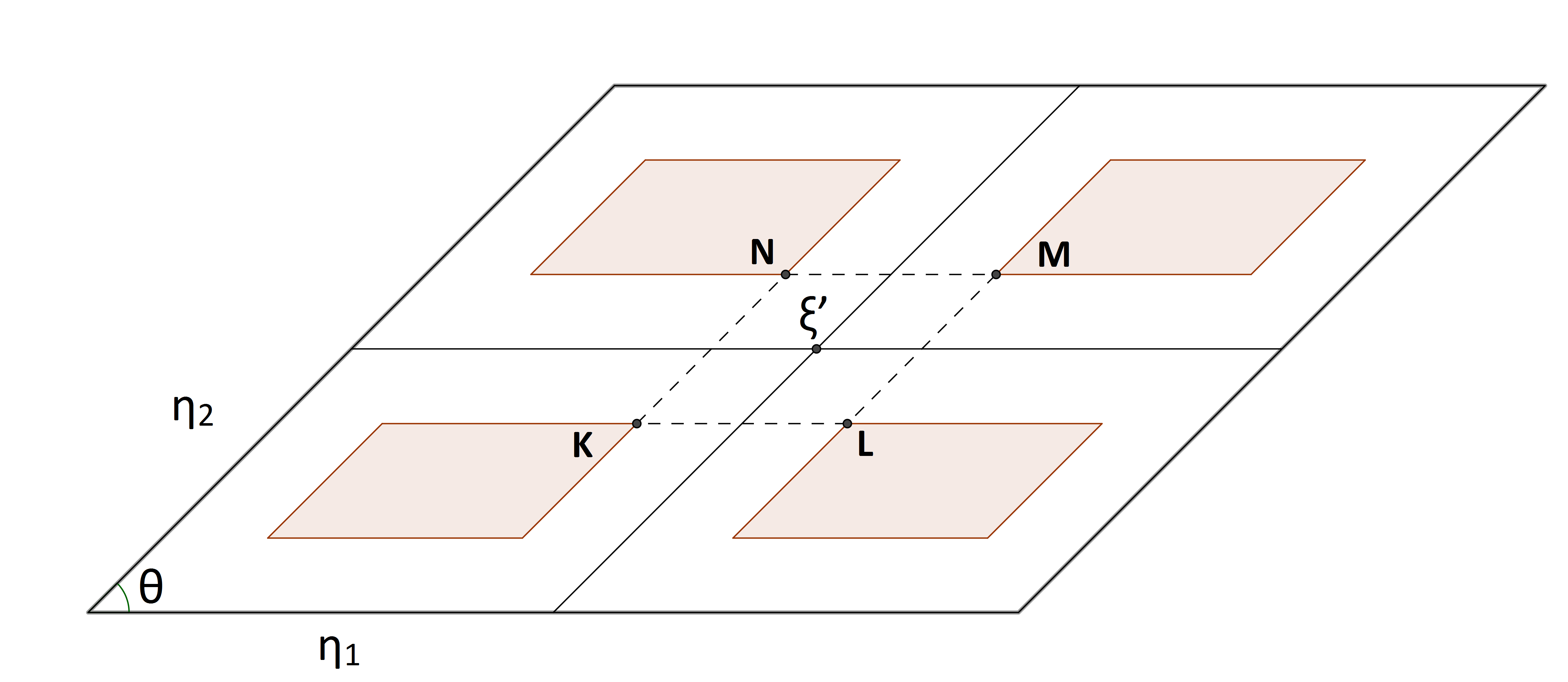}\label{ball1}
\end{figure}

Consider the region
$B_{-r}\left(\beta\zeta \X\right)$.
If $\beta w>2r$, 
then $B_{-r}\left(\beta\zeta \X\right) \neq \emptyset$.
Moreover, since this region $B_{-r}\left(\beta\zeta \X\right)$
has no intersection with any ball $B(x,r)$ centered at $x \in \C \setminus \beta\zeta \X$
of radius $r$,
the set $\left(z+\mathcal{L}\right)\cap \beta\zeta \X$ can not be empty and
gives an $r$-covering 
of $B_{-r}\left(\beta\zeta \X\right)$. 
That is, for each $z'\in B_{-r}\left(\beta\zeta \X\right)$,
there exists $d\in \L$ such that
$z+d \in \beta\zeta \X$ and the ball $B(z+d,r)$ contains $z'$.
As such, we see that $B_{-r/\beta}(\X)$
is $r/\beta$-covered by
$T^{-1}(z)$. 
Consequently, $B_{-r/\beta}(\X) +\L$ (see Figure \ref{ball1})  is $r/\beta$-covered by
$T^{-1}(z)+\L$.  
Now, we enlarge the radius $r/\beta$ to form 
a covering of the entire space $\C$.
To this end, we claim that extending the radius by a factor of
$\nu_2$ suffices.
From the inequality 
$\nu_2>1+1/\sin\theta$, we only have to check that a rhombus KLMN
in Figure \ref{ball1} determined by adjacent translates of 
$B_{-r/\beta}(\X)$
is covered.
Since $\nu_2$ is invariant under
$\theta \leftrightarrow \pi-\theta$, 
we prove the statement for $\theta \in (0,\pi/2]$.
Consider the Vorono\"i diagram of its four vertices $K,L,M$ and $N$. 
Then it can be seen easily that the 
minimum length required to achieve the goal is given by
the circumradius of the triangles $\Delta KLN$ and $\Delta LMN$, 
which are the acute triangles determined
by the smaller diagonal of the rhombus. This gives the constant 
$\nu_2$ and
proves the claim.
For an obtuse $\theta$, we have to switch to the other angle
$\pi-\theta$.
Refer to Figure \ref{ball2} below to compare the Vorono\"i diagrams
of two particular rhombuses.

\begin{figure}[htp]
\caption{Vorono\"i diagrams where $\alpha \in (0,\pi/2]$ and $\gamma \in (\pi/2,\pi)$}
	\centering
		\includegraphics[scale=0.05]{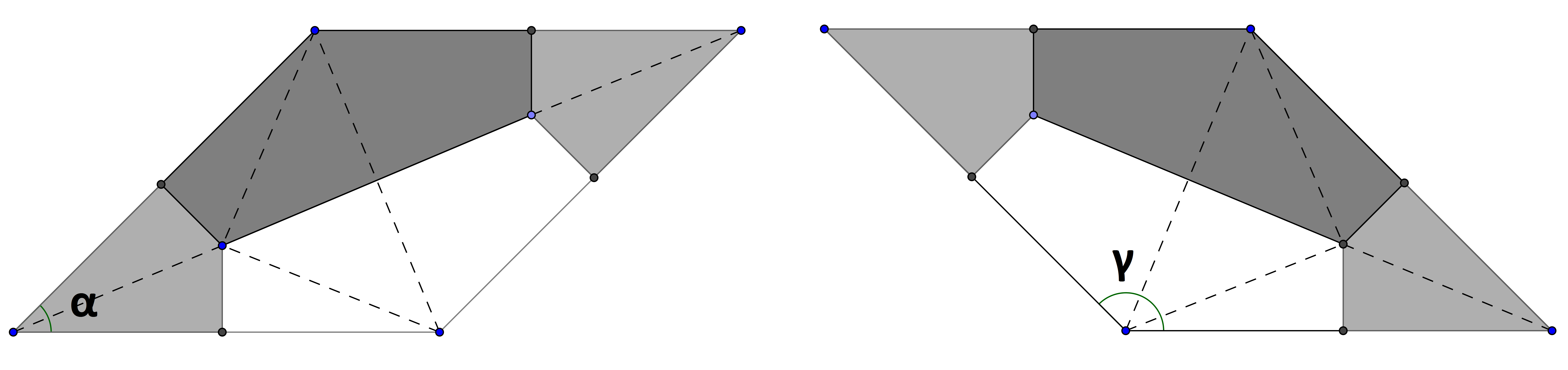}\label{ball2}
\end{figure}


Let $\beta > B_2$. We show by induction that for all $n\in \mathbb{N}$, 
$T^{-n}(z)$ provides an $r_n$-covering of $B_{-r_n}(\X)$, 
where 
$r_n=r\nu_2^{n-1}/\beta^n$.
Suppose this is the case for all $k\leq n$ for some $n\in \N$. 
We note that 
$T^{-(n+1)}(z)
=
\frac{1}{\beta \zeta} \left(\left(T^{-n}(z)+\L\right)\cap \beta\zeta \X\right)$. 
From $\beta > B_2$,
we have $r_n < r$.
Thus $\beta w> 2r_n $, implying that 
$B_{-r_n}\left(\beta\zeta \X \right) \neq \emptyset$. 
As $T^{-n}(z)+\L$ gives an $r_n$-covering of $B_{-r_n}(\X)+\L$, 
we can enlarge $r_n$ by a factor of $\nu_2$ to obtain a covering of $\C$,
and consequently, of $\beta \zeta \X$.
Now, for all $c\in \C\setminus \beta \zeta \X$, we have $B(c, r_n\nu_2) \cap B_{-r_n\nu_2}(\beta \zeta \X)=\emptyset$.
This implies that $(T^{-n}(z)+\L)\cap \beta \zeta \X$ is an $r_n\nu_2$-covering of $B_{-r_n\nu_2}(\beta \zeta \X)$. 
From this, it follows that $T^{-(n+1)}(z)$ is an $r_{n+1}$-covering of $B_{-r_{n+1}}(\X)$.
This finishes the induction
which completes the proof of the claim.

We continue to use the symmetry $\theta \leftrightarrow \pi-\theta$
and assume that $\theta\in (0,\pi/2]$.
In the course of the above proof, 
if we choose $z=\xi$, we can come up with a considerably finer covering of $\C$.
Observe that inside the parallelogram $KLMN$, there is a point $\xi'\in 
\xi+\L$ (see Figure \ref{ball1}).
The ball centered at $\xi'$ already covers a 
significant portion of the parallelogram. 
To proceed, we first note that
some rectangular strips along the perimeter of the translates of $\X$ can be covered
by balls $B(x,2r/\beta)$ where $x\in T^{-1}(\xi)+\L$
as shown in Figure \ref{ball7}.
Therefore, around $\xi'$, we need to cover a region comprising of four kite-shaped areas
given in Figure \ref{ball8}.
\begin{figure}[htp]
\caption{The rectangular strips}
	\centering
		\includegraphics[width=\textwidth]{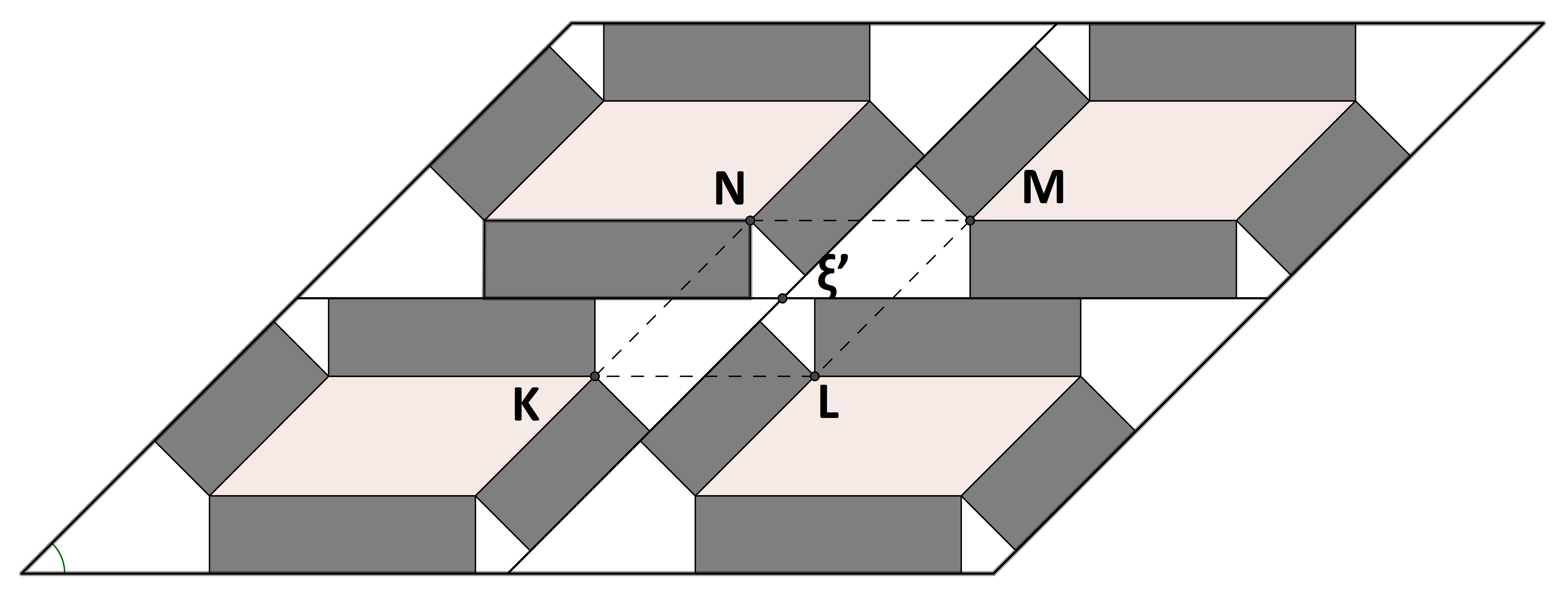}\label{ball7}
\end{figure}

\begin{figure}[htp]
\caption{The kites}
	\centering
		\includegraphics[width=\textwidth]{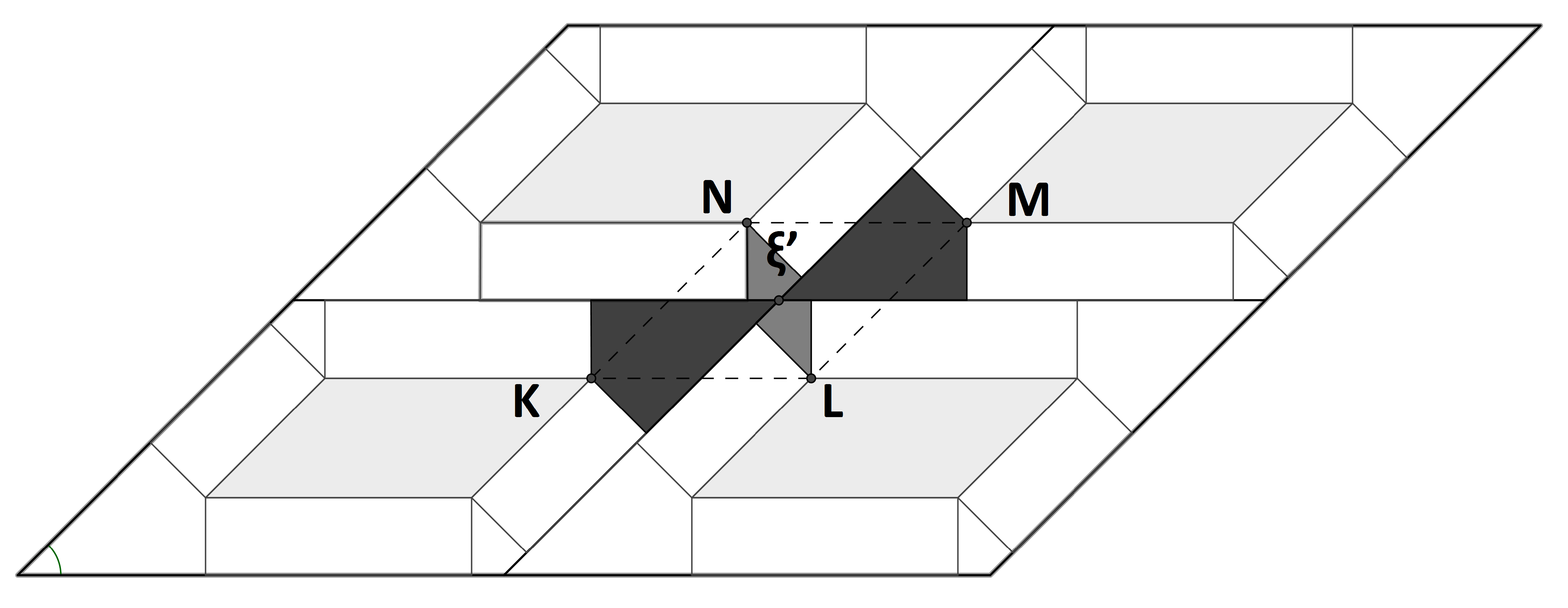}\label{ball8}
\end{figure}

Now, if $1/2<\tan (\theta/2)<2$, the ball $B\left(\xi',2r/\beta\right)$ contains
the bases of the two perpendiculars emanating from $M$
to the lines $\ell_1$ and  $\ell_2$, where
$\ell_j$ is the line parallel to $\eta_j$ 
and passing through $\xi'$ for $j=1,2$.
This means that the kite containing $M$ is covered by the balls
$B(\xi',2r/\beta)$ and $B(M,r/\beta)$. 
A similar argument shows that remaining kites are also covered.
Hence, we see that
$(\xi \cup T^{-1}(\xi))+\L$ gives a $2r/\beta$-covering of $\C$.

\begin{figure}[htp]
\caption{The balls $B(\xi',2r/\beta)$, $B(K,r/\beta)$ and $B(L,r/\beta)$}
	\centering
		\includegraphics[width=\textwidth]{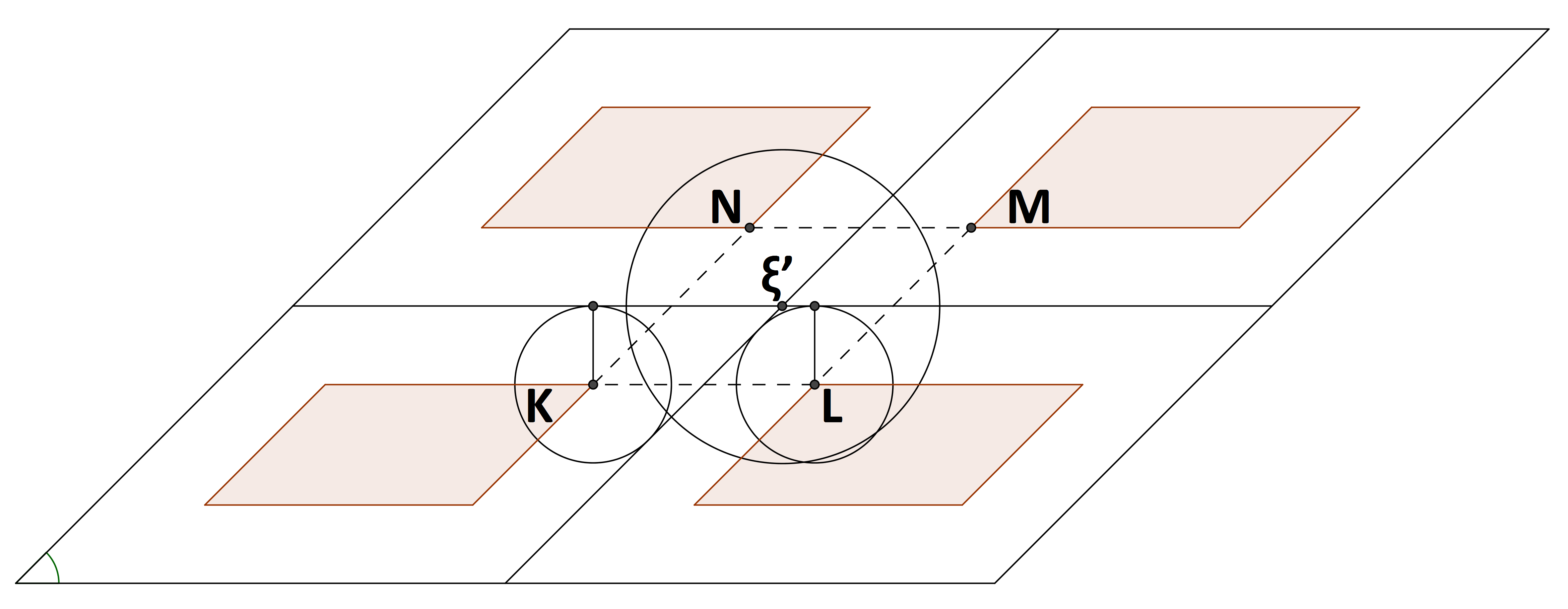}\label{ball3}
\end{figure}

For the other cases, we have to enlarge the radius a little more.
Figure \ref{ball3} shows such a case where there is a small remaining region yet to be covered. 
In Figure \ref{ball11}, 
we take a minimum $\rho>1$ 
such that the balls $B(\xi', (\rho +1)r/\beta)$ and $B(M,\rho r/\beta)$
intersect on the boundary of the kite.

A small computation yields that if 
$\rho=\frac{1+\cos \theta}{2(-1+\sin\theta+\cos \theta)}$,
$\C$ can be covered by balls centered at 
the elements of 
$(\xi \cup T^{-1}(\xi))+\L$ of radius $\rho r/ \beta$. 

\begin{figure}[htp]
\caption{Covering the kites}
	\centering

		\includegraphics[width=\textwidth]{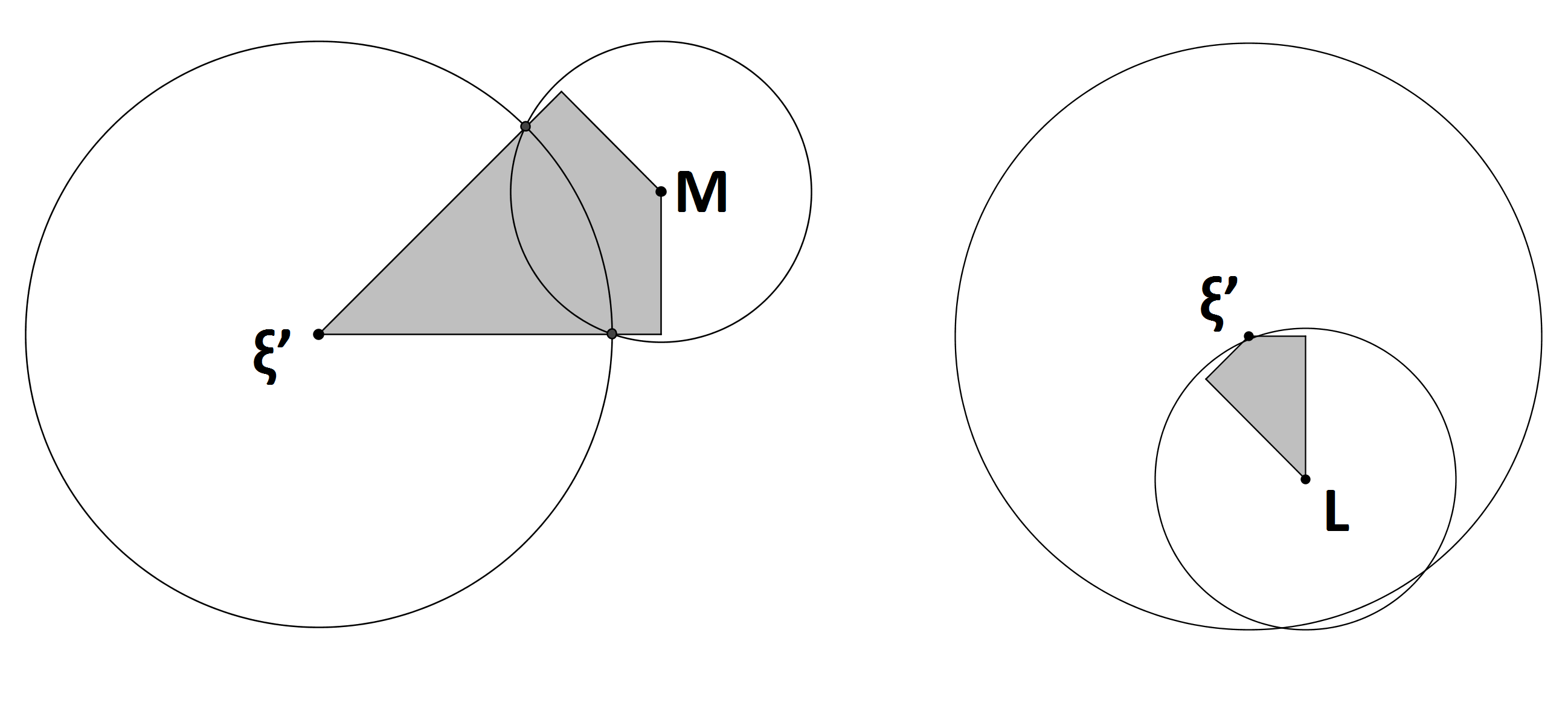}

\label{ball11}
\end{figure}

We can proceed with the same induction to see
that 
if $\beta > B_1$, $\bigcup_{j=0}^{n}T^{-j}(\xi)$ gives 
an $r \nu_1^{n-1}/\beta^n$-covering of $B_{-r \nu_1^{n-1}/\beta^n}(\X)$.


We saw that the choice $z=\xi$ makes the radius of the covering smaller. However
$\xi$ is unfortunately on the boundary of $\X$, which is not suitable for the later use.
So we select an appropriate $z\in \X$ which is very close to $\xi$. 
Since all inequalities in the above proof are open, one may find an
$\varepsilon_0>0$ that the first inductive step works 
for every point $z\in B(\xi,\varepsilon_0)$. 
Let $\Y:=\X \setminus\bigcup_{j=-\infty}^{\infty} T^j(\partial(\X))$
and select $z\in \Y$
such that 
$z \in B(\xi,\varepsilon_0 \nu_1^{n-2}/\beta^{n-1})$
for some integer $n\ge 2$.
This choice of $z$ is possible because
$\bigcup_{j=-\infty}^{\infty} T^j(\partial(\X))$
is a null set. 
Then the induction similarly works 
at least $n$ steps and we obtain the following statement.

\begin{it}
If $\beta>B_1$, then for  
any $\varepsilon>0$ there exists $z\in \Y$ and 
a positive integer $n$ such that $\bigcup_{j=0}^{n} T^{-j}(z)$ is an $\varepsilon$-covering of $B_{-\varepsilon}(\X)$.
\end{it}

We are ready to prove the first part of the theorem. Suppose $\beta > B_1$.
The proof of Theorem 5.2 in \cite{Saussol:00} implies that the support of
each ACIM contains an open ball where the associated Radon-Nikodym 
density has a positive lower bound.
Any such balls belonging to different ergodic ACIM's
must be disjoint. Let us 
assume that $\mu_j\ (j=1,2)$ are two different ergodic 
ACIM's of $(\X,T)$
with corresponding densities $h_j$.
Note that $h_j(z)>0$ implies $h_j(T(z))>0$ for almost all $z$ since
$h_j$ is a fixed point of the Perron Frobenius operator, whose associated 
Jacobian is positive and constant.
Choose open balls $B(x_j,s)$ such that $\einf_{B(x_j,s)}h_j>0$ for $j=1,2$.
From the above result, 
we can find $z\in  \Y$ and
positive integers $m_j$ $(j=1,2)$ such that
$T^{-m_j}(z)$ and $B(x_j,s)$ have a nontrivial intersection.
For $j=1,2$, let $u_j\in B(x_j,s)\cap T^{-m_j}(z)$. Then $T^{m_j}(u_j)=z$. 
Moreover, for some small balls $B(u_j,\delta_j)$ inside $B(x_j,s)$, 
we have 
\begin{eqnarray*}
T^{m_j}\left(B(u_j,\delta_j)\right)&=&B\left(T^{m_j}(u_j),\delta_j'\right)\\
																			&=&B(z,\delta_j'),
\end{eqnarray*}
where $\delta_j'>0$ is some small radius for $j=1,2$.
Therefore, $\einf_{B(z,\delta_j')} h_j>0$, which is a contradiction.
Thus, we see that the number of ergodic components is one, showing the 
first statement.

The second statement is subtler than the first one. Let $\beta >B_2$.
For $\varepsilon>0$, let
$$
N_{\varepsilon}= \{x\in \X\ |\ \einf_{B(x,\varepsilon)} h=0 \},
$$
where $h$ is the density of the ACIM $\mu$ and put
$N=\bigcap_{\varepsilon} N_{\varepsilon}$.
According to Proposition 5.1 in \cite{Saussol:00}, we know $\mu(N)=0$. 
We claim that $N$ is contained in 
$\bigcup_{j=-\infty}^{\infty} T^{j}(\partial(\X))$.
Assume 
that 
$z\not\in \bigcup_{j=-\infty}^{\infty} T^{j}(\partial(\X))$.
Choose $B(x,s)\subset \supp \mu$ with $x\in \X$ and $s>0$ such that
$\einf_{B(x,s)} h>0$.
Then there is a positive $n$ such that $T^{-n}(z) \cap B(x,s) \neq \emptyset$. This means
that there is a small ball $B(w,\varepsilon)\subset B(x,s)$ that $T^n(w)=z$ and
$T^n(B(w,\varepsilon))=B(z,\varepsilon')$.
However, $\einf_{B(z,\varepsilon')} h >0$ shows that
$z\not \in N$, which shows the claim. 
The claim implies $m(N)=0$
where $m$ is the 2-dimensional Lebesgue measure\footnote{
The proof of Proposition 5.1 in \cite{Saussol:00} guarantees $\mu(N)=0$. 
He wrote that this implies that $N$ is a null set (w.r.t. $\mu$), but 
it does not necessarily mean $m(N\cap \supp (\mu))=0$.
}.
Now we assume that $S\subset \X$ is measurable with 
$m(S)>0$. Since $m(S\setminus N)=m(S)>0$, take a Lebesgue density point $z\in 
S\setminus N$, i.e.,
$$
\lim_{t\rightarrow 0} \frac {m(B(z,t)\cap (S\setminus N))}{m(B(z,t))}=1.
$$
Since $z\not \in N$, there are positive $c$ and $\varepsilon_0$ such that
$\einf_{B(z,\varepsilon_0)} h >c$. Thus
$$
\mu(S)=\mu(S\setminus N)=\int_{S\setminus N} h dm > 
\int_{B(z,\varepsilon) \cap (S\setminus N) } c dm> 0,
$$
for a small $\varepsilon\le \varepsilon_0$ which shows that $m$ is absolutely continuous to $\mu$. 
\qed


\end{section}

\begin{section}{Proof of Lemma \ref{Sofic Definition}}\label{SD}
Recall that $\Y=\X\setminus \bigcup_{n=-\infty}^{\infty}
T^n(\partial(\X))$.
Define the set of {\it predecessors} associated to a point $z\in \Y$ by
$$
P(z)= \bigcup_{n=1}^{\infty}
\left\{ d(z')d(T (z'))\dots d(T^{n-1}(z')) \in \A^*
\ |\ z'\in T^{-n}(z) \right\},
$$
that is, the set of codings of all trajectories into $z\in \Y$
of the inverse images of the point $z$.
Introduce an equivalence relation $z_1 \sim z_2$ by $P(z_1)=P(z_2)$.
It is clear that the cardinality of equivalence classes in $\Y/\sim$ 
is finite if and only if
the system is sofic (cf. \cite[Theorem 3.2.10]{Lind-Marcus:95}).
By the definition of the map $T$, it 
is plain to see by induction on $K$, 
that $\X\setminus \bigcup_{n=1}^{K} T^n(\partial(\X))$ consists of finite number of
open polygons and
each end point of a discontinuity segment 
must be on {\it another segment of a different slope}\footnote{If a segment of $T^n(\partial(\X))$ falls into $\partial(\X)$,
then we discard the segment, because the soficness 
is defined over $\Y$.}.
An open polygon may be cut into two or more
pieces by a broken line of $T^{K+1}(\partial(\X))$.
We see that any points $x$ and $y$ separated by the broken line 
are inequivalent, as one of $P(x)$ and $P(y)$
has at least one more predecessor than the other.
Suppose that $\bigcup_{n=1}^{\infty} T^n(\partial(\X))$ 
is an infinite union of segments. Then as we increase $K$ by $1$,
at least one open polygon of
$\X\setminus \bigcup_{n=1}^{K} T^n(\partial(\X))$
is separated by a broken line coming from
$T^{K+1}(\partial(\X))$. 
In fact, if not then $T^{K+1}(\partial(\X))$ must be totally 
contained in $Q:=\partial(\X)\cup \bigcup_{n=1}^{K} T^n(\partial(\X))$ 
and we have $T^{m}(\partial(\X))\subset Q$ for $m\ge K+1$. However 
there are only finitely many segments whose end points
lie on other segments of different slopes in $Q$, which shows that the sequence
$(T^{m}(\partial(\X)))\ (m>K)$ is eventually periodic,
giving a contradiction. 
Consequently we always find an additional equivalent class through $K\rightarrow K+1$. 
This shows that the system can not be sofic.

For the reverse implication, we consider the
partition of $\X$ into finitely many disjoint polygons 
induced by $\bigcup_{n=0}^{\infty}T^n(\partial(\X))$. 
Taking discontinuities into account, 
such polygons may not be open nor closed. 
Let $P_1, \ldots, P_r$ be the polygons in the partition.
It is clear that for $j\in \{1, \ldots,r\}$, 
$T(P_j)=\bigcup_{k\in \mathcal{I}}P_k$
for some $\mathcal{I}\subseteq \{1, \ldots, r\}$. This follows
from the fact that the set $\bigcup_{n=0}^{\infty} T^n(\partial(\X))$
is $T$-invariant. 
For $d \in \mathcal{A}$, let 
$$[d]:=\{z\in \X| d_1(z)=d\}.$$ 
Suppose $P_j \cap [d] \neq \emptyset$. 
Since $\beta \zeta [d]=\beta \zeta \X \cap (\X+d)$, then the boundary of $T([d])$ lies in $\bigcup_{n=0}^{1}T^n(\partial(\X))$.
Note that 
$$T(P_j \cap [d]) \subseteq T(P_j) = \bigcup_{k\in \mathcal{I}}P_k$$ 
and
\begin{eqnarray*}
\beta \zeta (P_j \cap [d]) &=& \beta \zeta P_j \cap \beta \zeta [d]\\
													 &=& \beta \zeta P_j \cap (\X+d).
\end{eqnarray*}
Thus, $T(P_j \cap [d])=\bigcup_{k\in \mathcal{I}^*}P_k$ where $\mathcal{I}^* \subseteq \mathcal{I}$.
From the partition, we define a labeled directed graph $G$.
Let $$V(G):=\{P_1, \ldots, P_r\}$$ be the
vertex set of $G$. We build the edge set and define the labeling as follows. 
For $j,k\in \{1,\ldots, r\}$ and $d\in \mathcal{A}$,
there is an edge labeled $d$ from $P_j$ to $P_k$ if 
$P_k$ is contained in $T(P_j \cap [d])$. It is clear that $G$ is a sofic graph 
describing $(\X,T)$.
\qed

\begin{rem}  The sofic shift obtained in the latter part of the above proof
is irreducible if $(\X,T)$ admits the ACIM equivalent to the Lebesgue measure.
By construction, the resulting labeled graph is the 
minimum left resolving presentation of the irreducible sofic shift.
Therefore it is easy to check whether the system is a shift of finite type
or not by checking synchronizing words through 
backward reading of the graph (see \cite[Theorem 3.4.17]{Lind-Marcus:95}). 
\end{rem}

\end{section}

\begin{section}{Proof of Theorem \ref{Sofic}}\label{S}
We have to study the growth of
$\bigcup_{n=1}^{K} U^n\left(\partial\left(\left[0,1\right)^2\right)\right)$
as $K$ increases. Our idea is to record only the information 
of the set of {\it lines} which include this finite union of segments. 
Thus, we are interested in studying the union of the lines containing the
segments
whose defining equations are 
of the form $f(X,Y)=(A,B) \begin{pmatrix} X\\Y \end{pmatrix}+C$, where 
$(0,0)\neq (A,B)\in \mathbb{R}^2$. 
We often identify the line and its defining equation.  
Then the image under $U$ of the line 
is given by the defining equation 
$f(X', Y')=0$ with
$$
\begin{pmatrix}X\\Y\end{pmatrix}
= \beta \begin{pmatrix}a_{11}& a_{12} \\ a_{21} & a_{22} \end{pmatrix}
\begin{pmatrix} X' \\ Y'\end{pmatrix} - \begin{pmatrix} c_1 \\ c_2\end{pmatrix}
$$
where
$$
\begin{pmatrix} c_1 \\ c_2\end{pmatrix} \in 
\Delta:=\left. \left\{
\begin{pmatrix}\lfloor \beta (a_{11}x+a_{12}y)+b_1 \rfloor-b_1\\
\lfloor \beta (a_{21} x + a_{22}y)+b_2 \rfloor-b_2\end{pmatrix}
\right| 
 0\leq x,y< 1
\right\}.
$$
Since $\Delta$ is a bounded set of lattice points, it is 
a finite set.
As multiplication by $\zeta$ acts as $q$-fold rotation on $\C$, we have
\begin{equation}
\label{RotRep}
\begin{pmatrix}a_{11}& a_{12} \\ a_{21} & a_{22} \end{pmatrix}^{-1}
=\begin{pmatrix}a_{22}& -a_{12} \\ -a_{21} & a_{11}\end{pmatrix},
\quad
\begin{pmatrix}a_{11}& a_{12} \\ a_{21} & a_{22} \end{pmatrix}^{q}
=
\begin{pmatrix}1&0 \\ 0 & 1 \end{pmatrix}.
\end{equation}
Therefore the image of the line under $U$ is
$$
\frac 1\beta (A,B) 
\begin{pmatrix}a_{22}& -a_{12} \\ -a_{21} & a_{11}\end{pmatrix} 
\begin{pmatrix} X+c_1 \\ Y+c_2\end{pmatrix} + C =0.
$$
Multiplying by $\beta$, we obtain a correspondence of the 
coefficient vectors of the defining equations:
\begin{equation}\label{Corr}
\left(A^{(n)},B^{(n)},C^{(n)}\right) \rightarrow \left(A^{(n+1)},B^{(n+1)},C^{(n+1)}\right)
\end{equation}
where
\begin{equation}
\label{NextAB}
\left(A^{(n+1)},B^{(n+1)}\right)=\left(A^{(n)},B^{(n)}\right)
\begin{pmatrix} a_{22}& -a_{12} \\ -a_{21} & a_{11}\end{pmatrix},
\end{equation}
\begin{equation}
\label{NextC}
C^{(n+1)}= \beta C^{(n)} + \left(A^{(n)},B^{(n)}\right) 
\begin{pmatrix} a_{22}& -a_{12} \\ -a_{21} & a_{11}\end{pmatrix} 
\begin{pmatrix} c_1 \\ c_2\end{pmatrix}
\end{equation}
with $\left(A^{(0)},B^{(0)},C^{(0)}\right)=(A,B,C)$. 
Note that (\ref{Corr}) is not one-to-one, since we have many 
choices for $\begin{pmatrix} c_1 \\ c_2\end{pmatrix}$ from $\Delta$. 
Here we introduce an obvious restriction on $C^{(n)}$ that four values
$$
\left.\left\{ A^{(n)} s+ B^{(n)}t +C^{(n)} \ \right| (s,t)\in \{(0,0),(1,0),(0,1),(1,1)\} 
\right\} 
$$
are not simultaneously positive nor negative, to ensure
that the resulting lines intersect the closure of $\X$.
All the same we have to note that
the resulting lines may contain irrelevant ones\footnote{
Therefore the resulting lines are {\it potential} discontinuities.
In the actual algorithm to obtain the associated graph of the sofic shift, 
it is simpler to abandon such irrelevant lines at each step.
However in doing so, 
we have to record the position of end points of discontinuity segments, 
which makes the process involved.
}
which do not actually contain 
a segment of
$\bigcup_{j=0}^{\infty} U^j(\partial([0,1)^2))$.
From (\ref{RotRep}), $(A^{(n)}, B^{(n)})$ is clearly periodic
with period $q$, and our task is to prove that the set of all 
$C^{(n)}$ given by this iteration is finite.
We call the set $\overline{U} \subset \mathbb{Q}(\beta)$ of all the $C^{(n)}$'s 
arising from $\partial([0,1)^2)$, together with 0 and -1, the set of {\it intercepts} 
of $U$. 

Let $\beta_1=\beta, \beta_2,\ldots, \beta_{d}$ be the conjugates of $\beta$. 
For $k=1,\dots,d$, define $\sigma_{k}:\mathbb{Q}(\beta)\rightarrow\mathbb{Q}(\beta_k)$
to be the conjugate map that sends $\beta$ to $\beta_{k}$. 
To demonstrate the finiteness of $\overline{U}$, 
we show that $\sigma_k\left(C^{(n)}\right)$ 
is bounded for $k=1,\ldots,d$. 
From (\ref{NextAB}),(\ref{NextC}) and (\ref{RotRep}), we have 
$
C^{(n+1)}= \beta C^{(n)} + m
$
where $m$ is an element of 
$$
M:= \left\{ \left. (A,B)
\begin{pmatrix} a_{22}& -a_{12} \\ -a_{21} & a_{11}\end{pmatrix}^n 
\begin{pmatrix} c_1 \\ c_2\end{pmatrix} \ \right| 
(A,B)\in \{ (0,1),(1,0) \}, 
n=0,1,\dots, q-1,  \begin{pmatrix} c_1\\c_2\end{pmatrix} \in \Delta\right\}.
$$
Here we use the fact that $(A,B,C)\in \{(1,0,0),(1,0,-1),(0,1,0),(0,1,-1)\}$
gives $\partial([0,1)^2)$.
By the finiteness of $\Delta$, $M \subset \Q(\beta)$ 
is also a finite set. 
Taking a common denominator, there is a fixed $N\in \N$ that
$C^{(n)} \in \frac{1}{N}\Z [\beta]$.
Let 
$\omega_k := \max \{1, \max\limits_{m\in M} \{|\sigma_k(m)|\}\}$. 
Then, if $k=2,\ldots, d$, we have 
$$
\left|\sigma_{k}\left(C^{(n)}\right)\right| 
\leq |(\beta_k)^n| + \omega_k \sum_{j=0}^{n-1}|\beta_k|^j 
\leq \frac{\omega_k}{1-|\beta_k|}.
$$
For $k=1$, since the line 
$A^{(n)}X+B^{(n)}Y+C^{(n)}=0$ passes through
$[0,1]^2$,
it follows that 
$$
\left|\sigma_1\left(C^{(n)}\right)\right|=
\left|C^{(n)}\right| \leq \max_{l=0,1,\dots q-1} \left(\left|A^{(l)}\right|+\left|B^{(l)}\right|\right).
$$
by the periodicity of $A^{(n)}$ and $B^{(n)}$.
\qed

\end{section}

\begin{section}{Proof of Theorem \ref{Non-Sofic}}\label{Non}
Put $\omega=(1+\sqrt{5})/2$. 
From $\xi=0$, $\eta_1=1$, $\eta_2=\zeta=\exp(2\pi i/5)$
and a trivial relation $\zeta^2=(\zeta+\zeta^{-1})\zeta-1$,
we have $b_1=b_2=0$ and $a_{11}=0, a_{12}=-1, a_{21}=1, a_{22}=1/\omega$.
Therefore, we have
$$
U\left(\begin{pmatrix}x\\y\end{pmatrix}\right)= 
\begin{pmatrix} -\beta y-\lfloor -\beta y \rfloor\\ 
\beta (x+y/\omega)  
-\lfloor \beta (x+y/\omega) \rfloor
\end{pmatrix}.
$$
Clearly, $\sqrt{5}\not \in \Q(\beta)$ is equivalent to
 $\Q(\beta)\cap \Q(\omega)=\Q$.
Since $\Q(\omega)$ is a Galois extension over $\Q$,
this implies that $\Q(\omega)$ and $\Q(\beta)$
are linearly disjoint and there exists a conjugate map 
$\sigma\in \mathrm{Gal}(\Q(\beta,\omega)/\Q(\beta))$ with
$\sigma(\beta)=\beta$ and $\sigma(\omega)=-1/\omega$.

From (\ref{NextAB}) and (\ref{NextC}) we see, 
$$
C^{(n+1)}=\beta C^{(n)}+\left(Ac_{11}^{(n+1)}+Bc_{21}^{(n+1)}\right)c_1+
\left(Ac_{12}^{(n+1)}+Bc_{22}^{(n+1)}\right)c_2\in \overline{U}$$ 
for some $\left( \begin{matrix} c_1\\ c_2\\\end{matrix}\right)\in \Delta$
and
$$
\begin{pmatrix} c_{11}^{(n)} & c_{12}^{(n)}\\ c_{21}^{(n)} & c_{22}^{(n)} 
\end{pmatrix}
= \begin{pmatrix} a_{22} & -a_{12}\\ -a_{21} & a_{11}\end{pmatrix}^n.
$$
Consider the case where $(A,B,C)=(1,0,-1)$. 
Then, $$C^{(n+1)}=\beta C^{(n)}+c_{11}^{(n+1)}c_1+c_{12}^{(n+1)}c_2.$$
Applying $\sigma$, we get 
$\sigma \left(C^{(n+1)}\right)=\beta \sigma\left(C^{(n)}\right)+\sigma\left(c_{11}^{(n+1)}c_1+c_{12}^{(n+1)}c_2\right)$.
It follows that
\begin{eqnarray*}
\left|\sigma \left(C^{(n+1)}\right)\right| &\geq&
\beta\left| \sigma\left(C^{(n)}\right)\right|
-\left|\sigma\left(c_{11}^{(n+1)}c_1+c_{12}^{(n+1)}c_2\right)\right|\\
&=&\beta\left| \sigma\left(C^{(n)}\right)\right|
-\left|\sigma\left(c_{11}^{(n+1)}\right)c_1+\sigma\left(c_{12}^{(n+1)}\right)c_2\right|\\
&\geq&\beta\left| \sigma\left(C^{(n)}\right)\right|-D,
\end{eqnarray*}
where 
\begin{eqnarray*}
D&:=&\max \limits_{n\in \N} 
\max \limits_{\Delta}
\left\{ \left| \sigma \left(c_{11}^{(n)}\right) c_1 + \sigma \left(c_{12}^{(n)}\right) c_2 \right| \right\}\\
&\leq&\max \limits_{n\in \N} 
\max \limits_{\Delta}
\left\{\left|\sigma\left(c_{11}^{(n)}\right)\right|\left|c_1\right|+\left|\sigma\left(c_{12}^{(n)}\right)\right|\left|c_2\right|\right\}.
\end{eqnarray*}
Direct computation yields
$$\left(\sigma\left(c_{11}^{(n)}\right),\sigma\left(c_{12}^{(n)}\right)\right)=
\begin{cases}
(1,0) & n\equiv 0\pmod{5}\\
(-\omega,1) & n\equiv 1\pmod{5}\\
(\omega,-\omega) & n\equiv 2\pmod{5}\\
(-1,\omega) & n\equiv 3\pmod{5}\\
(0,-1) & n\equiv 4\pmod{5}.
\end{cases}$$ 
Hence, $D\leq \omega \max \limits_{\Delta}
\left\{\left|c_1\right|+\left|c_2\right|\right\}=\omega \left(\left\lfloor \beta \omega\right\rfloor+\left\lceil \beta\right\rceil\right)$.
Accordingly, for all $n \in \N$,
$$\left|\sigma \left(C^{(n+1)}\right)\right| \geq
\beta\left| \sigma\left(C^{(n)}\right)\right|-\omega \left(\left\lfloor \beta \omega\right\rfloor+\left\lceil \beta\right\rceil\right).$$
Therefore, if $$\left|\sigma \left(C^{(n)}\right)\right|>\frac{\omega \left(\left\lfloor \beta \omega\right\rfloor+\left\lceil \beta\right\rceil\right)}{\beta - 1}$$
for some $n \in \N$, then $\left\{\sigma \left(C^{(n)}\right)\left. \right|n\in \N \right\}$ diverges.
Now, it is easy to check that $\left(A^{(1)},B^{(1)},C^{(1)}\right)=(
\omega -1, 1, (\omega -1)\left\lfloor 
-\beta \right\rfloor+\left\lfloor \beta \omega\right\rfloor-\beta)$ 
gives a line which actually includes a discontinuity segment.
Under the assumption $\beta >\frac{1}{4} \left(13+3 \sqrt{5}-\sqrt{70-2 \sqrt{5}}\right)
\approx 2.90332$, we have
\begin{eqnarray*}
\left|\sigma \left((\omega -1)\left\lfloor -\beta \right\rfloor+\left\lfloor \beta \omega\right\rfloor-\beta\right)\right|
&=&\sigma \left((\omega -1)\left\lfloor -\beta \right\rfloor+\left\lfloor \beta \omega\right\rfloor-\beta\right)\\
&=&-\omega \left\lfloor -\beta \right\rfloor+\left\lfloor \beta \omega\right\rfloor-\beta,\end{eqnarray*} 
and
$$-\omega \left\lfloor -\beta \right\rfloor+\left\lfloor \beta \omega\right\rfloor-\beta > \frac{\omega \left(\left\lfloor \beta \omega\right\rfloor+\left\lceil \beta\right\rceil\right)}{\beta - 1}.$$ We therefore conclude that $\left\{\sigma \left(C^{(n)}\right)\left.\right|n\in\N\right\}$ is unbounded.
Now we have shown that once we had chosen $C^{(1)}$ as above, 
for {\it every} possible sequence $\left(C^{(n)}\right)$, its conjugate sequence 
$\left(\sigma\left(C^{(n)}\right)\right)\ (n=1,2,3,\dots)$ diverges. This implies that the set 
of discontinuities can not be finite.
\qed

\end{section}

\begin{section}{Examples}

Taking $\beta$ small, we can find a family of systems $(\X,T)$ with more than one ACIM.

\begin{exam}
\label{NonErgode}
Let $\zeta=i$ and $\beta = 1.039$. Set $\eta_1=2.92$,
$\eta_2=\exp(\pi i/3)$ and $\xi = 0$. The 
distribution of eventual orbits of $T$ of randomly chosen points  
is depicted in Figures \ref{TwoCompo1} and \ref{TwoCompo2}. From these figures, it is not difficult to make explicit the polygons
bounded by horizontal and vertical segments (easier after filling holes)
within which restrictions of $T$ are well-defined. 

\begin{figure}[htp]
\centering
\includegraphics[scale=1]{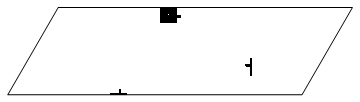}	
\caption{First Component\label{TwoCompo1}}
\end{figure}

\begin{figure}[htp]
\centering
\includegraphics[scale=1]{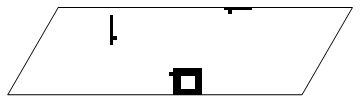}
\caption{Second Component\label{TwoCompo2}}
\end{figure}

\noindent
This leads us to 
a rigorous proof of the existence of two distinct ergodic components.
In Figure \ref{TwoCompo1}, the largest polygon is composed of two shapes $E$ and $F$
as in Figure \ref{EF}. The ratio of two sides of the rectangle $E$ is $1:\beta$.
By successive applications of $T$, the four vertices of $E$ are easily computed:
$$
x+\frac{\sqrt{3}i}2,\ x+yi,\ 
x + \frac 1{\beta} \left(\frac{\sqrt{3}}2-y\right)+yi,\ 
x + \frac 1{\beta}\left(\frac{\sqrt{3}}2 - y\right)+\frac{\sqrt{3}i}2
$$
with
$x=\eta_1-\frac{\sqrt{3}}{2}\beta -\frac 12$
and $y=\beta x-\frac {\sqrt{3}}2$ and the vertices of $F$ 
in 
counter-clockwise ordering 
are
$$
x + \frac 1{\beta} \left(\frac{\sqrt{3}}2-y\right)+yi,\
\gamma+yi,\ 
\gamma+vi,\ 
u+vi,
$$
$$
u+v'i,\ 
\gamma+v'i,\  
\gamma+\frac{\sqrt{3}i}2,\ 
x + \frac 1{\beta}\left(\frac{\sqrt{3}}2 - y\right)+\frac{\sqrt{3}i}2
$$
with $\gamma=-\beta y+\eta_1-\frac 12$, 
$u+vi=T^3(\gamma+\sqrt{3}i/2)$ 
and
$u+v'i=T^2(x-1/2)$.
Two other polygons found in Figure \ref{TwoCompo1}
are $T(F)$ and $T^2(F)$, which are similar to $F$ with the ratio $\beta$ and $\beta^2$.
We readily confirm the set equation $E \cup F = T(E) \cup T^3(F)$ (see Figure \ref{SE}). 
Hence the restriction of $T$ is well defined on the set
$$
Y:=E \cup F \cup T(F) \cup T^2(F)
$$
and defines a piecewise expanding map. Thus there is at least one ACIM whose
support is contained in $Y$. The same discussion 
can be done for Figure \ref{TwoCompo2}. 
The resulting supports of ACIM's are clearly disjoint.

\begin{figure}[htp]
\centering
\includegraphics[scale=0.6]{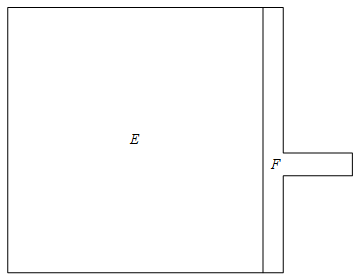}	
\caption{$E$ and $F$ \label{EF}}
\end{figure}

\begin{figure}[htp]
\centering
\includegraphics[scale=0.6]{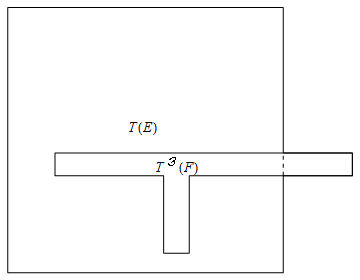}	
\caption{Confirmation of the set equation \label{SE}}
\end{figure}

The same situation happens when $\beta$ and $\eta_1$ satisfy
$$
\frac {\sqrt{3}}2 \beta+1+\frac {\sqrt{3}}{\beta}-\frac {\sqrt{3}}{2\beta^3}
\leq \eta_1 \leq \frac 12+ \frac {\sqrt{3}}{\beta}+\frac {\sqrt{3}}{2\beta^3}
$$
while other parameters are fixed. 
The corresponding region is shaded in Figure \ref{Param}.
This example gives an uncountable family of systems
with at least two ACIM's. 

\begin{figure}[htp]
\centering
\includegraphics[scale=0.7]{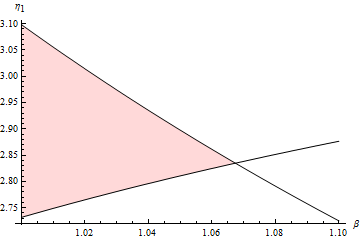}	
\caption{Non ergodic parameters \label{Param}}
\end{figure}

\end{exam}

In the following we give some examples of sofic systems.

\begin{exam}
Let $\zeta=\exp(2\pi i/3)$ and $\beta = 1+\sqrt{2}$. Set $\eta_1=1$ and $\eta_2=\zeta ^2$ and $(\beta\zeta -1)\xi = 3-\beta$. 
From $r(\L)=1/\sqrt{3}$, $w(\X)=\sqrt{3}/2$ we have $\beta>B_2=7/3$ and
there is a unique ACIM equivalent to Lebesgue measure by Theorem \ref{Ergodic}.
We consider the symbolic dynamical system associated to the map $T$. 
The set $\mathcal{A}$ is given by

\begin{align*}
&\bigl\{a=-1-\zeta^2,
b=-\zeta^2,
c=1-\zeta^2,
d=2-\zeta^2,
e=-2-2\zeta^2,
f=-1-2\zeta^2,\\
&g=-2\zeta^2,
h=1-2\zeta^2,
j=-2-3\zeta^2,
k=-1-3\zeta^2,
l=-3\zeta^2 \bigr\}.
\end{align*}

\begin{figure}[htp]
\centering
\includegraphics[scale=0.08]{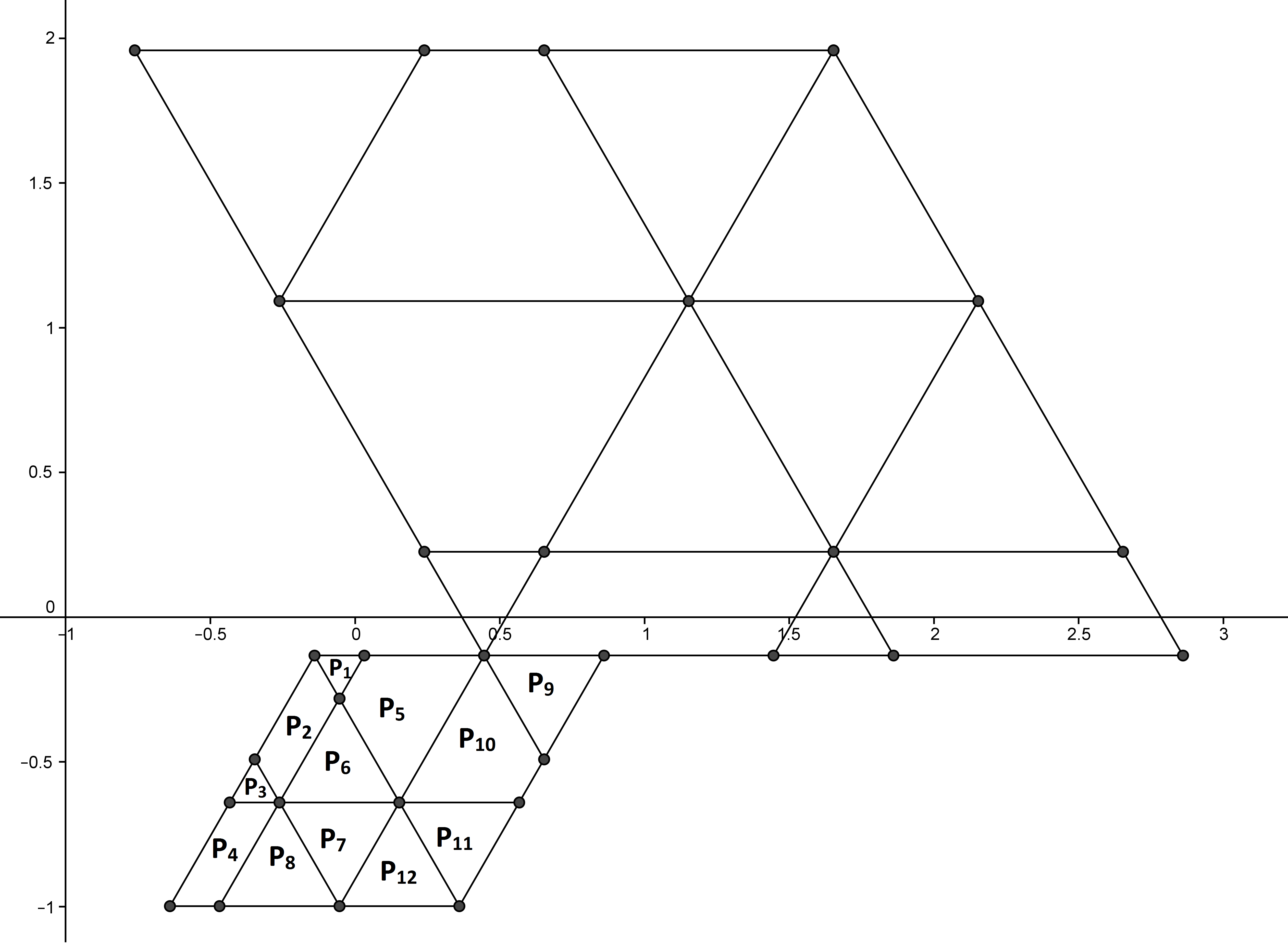}	
\caption{$\X$ and $\beta \zeta \X$\label{fig:threefold}}
\end{figure}

In Figure \ref{fig:threefold}, we see that the discontinuity
lines are finite and partition the fundamental domain $\X$ into
disjoint components $P_n$, $n=1,\ldots,12$. We also see in the figure
the expanded fundamental region $\beta \zeta \X $.

It is easy to confirm
that the image of $P_n$ under $T$ is given by Table \ref{three}.
From this table, we construct the sofic graph (see Figure \ref{fig:soficgraph}) 
as described in \S\ref{SD}.

\begin{table}[ht]
\caption{} 
\centering 
\begin{tabular}{l l c} 
\hline
\hline 
$n$ &  $T(P_n)$ & $\delta$  \\ [0.5ex]  

\hline 
1 &  $P_{7}$ & $b$ \\ 
2 &  $P_{11} \cup P_{12}$ & $b$  \\
  &  $P_{4} \cup P_{7}\cup P_8$ & $c$\\
3 &  $P_{12}$  &   $c$ \\
4 &  $P_{11}$ &   $c$ \\
	&  $P_{4} \cup P_{7}\cup P_8\cup P_{12}$ & $d$\\
5 &  $P_9$ &   $a$\\
	&  $P_{1} \cup P_{2}\cup P_5 \cup P_6$ & $b$ \\
	&  $P_{7} \cup P_{11}\cup P_{12}$ & $f$  \\
	&  $P_{4} \cup P_{7}\cup P_8$ & $g$  \\
6 &  $P_{9}\cup P_{10}$ & $b$ \\
	&	 $P_2\cup P_3 \cup P_6$ & $c$\\
	&  $P_{12}$ & $g$  \\
7 &  $P_{1} \cup P_{5}$ & $c$\\
  &  $P_{11}$ & $g$\\
  &  $P_{4} \cup P_{7}\cup P_8$ & $h$  \\
[1ex]\hline 
\end{tabular}
\qquad
\begin{tabular}{l l c} 
\hline
\hline 
$n$ &  $T(P_n)$ & $\delta$  \\ [0.5ex]  

\hline 
8 &  $P_{9}\cup P_{10}$ & $c$ \\
	&  $P_{2} \cup P_{3}\cup P_6$ & $d$  \\
	&  $P_{12}$ & $h$  \\
9 &  $P_{9}$ & $e$\\
	&  $P_{1} \cup P_{2}$ & $f$  \\
	&  $P_{7} \cup P_{11}\cup P_{12}$ & $j$  \\
	&  $P_{4}$ & $k$  \\
10  &  $P_{5}\cup P_{6}\cup P_{9}\cup P_{10}$ & $f$ \\
    &  $P_{2}\cup P_{3}\cup P_{6}$ & $g$ \\
		&  $ P_{7}\cup P_8\cup P_{12}$ & $k$  \\
11	&  $P_{1}\cup P_{5}$ & $g$ \\
		&  $P_{11}$ & $k$ \\
		&  $P_{4} \cup P_{7}\cup P_8$ & $l$  \\
12	&  $P_{9} \cup P_{10}$ & $g$  \\
		&  $P_{2} \cup P_{3}\cup P_6$ & $h$  \\
		&  $P_{12}$ & $l$  \\
[1ex]\hline 
\end{tabular}
\label{three}
\end{table}

\begin{figure}
\centering
\begin{tikzpicture}[midarrow/.style={decoration={markings, mark=at position 0.5 with {\arrow{stealth}}},postaction={decorate}}]
    \node[draw,circle] (r1) at (-2.5,-2.5) {$P_1$};
    \node[draw,circle] (r2) at (7.5,-2.5) {$P_2$};
    \node[draw,circle] (r3) at (10,-2.5) {$P_3$};
    \node[draw,circle] (r4) at (2.5,0) {$P_4$};
		\node[draw,circle] (r5) at (0,-2.5) {$P_5$};				
		\node[draw,circle] (r6) at (2.5,-5) {$P_6$};
		\node[draw,circle] (r7) at (0,0) {$P_7$};
		\node[draw,circle] (r8) at (5,0) {$P_8$};			
		\node[draw,circle] (r9) at (0,-5) {$P_9$};
		\node[draw,circle] (r10) at (5,-5) {$P_{10}$};
		\node[draw,circle] (r11) at (2.5,-2.5) {$P_{11}$};
		\node[draw,circle] (r12) at (5,-2.5) {$P_{12}$};

\draw[midarrow,bend right=15] (r1) to node[below] {$$} (r7);

\draw[midarrow,bend right=15] (r2) to node[above] {$$} (r11);
\draw[midarrow,bend right=15] (r2) to node[above] {$$} (r12);
\draw[midarrow,bend right=15] (r2) to node[above] {$$} (r4);
\draw[midarrow,bend right=15] (r2) to node[above] {$$} (r7);
\draw[midarrow,bend right=15] (r2) to node[above] {$$} (r8);

\draw[midarrow,bend right=15] (r3) to node[above] {$$} (r12);

\draw[midarrow,bend right=15] (r4) to node[above] {$$} (r11);		 	
\draw[midarrow,looseness=15] (r4) to[out=110,in=150] node[left] {$$} (r4);		 	
\draw[midarrow,bend right=15] (r4) to node[above] {$$} (r7);		 	
\draw[midarrow,bend left=15] (r4) to node[below] {$$} (r8);		
\draw[midarrow,bend right=15] (r4) to node[above] {$$} (r12);

\draw[midarrow,bend right=15] (r5) to node[above] {$$} (r9); 	
\draw[midarrow,bend left=15] (r5) to node[above] {$$} (r1); 	
\draw[midarrow,bend right=15] (r5) to node[above] {$$} (r2); 	
\draw[midarrow,looseness=15] (r5) to[out=110,in=150] node[left] {$$} (r5); 	
\draw[midarrow,bend right=15] (r5) to node[above] {$$} (r6); 	
\draw[midarrow,bend right=15] (r5) to node[above] {$$} (r7); 	
\draw[midarrow,bend right=15] (r5) to node[above] {$$} (r11); 	
\draw[midarrow,bend right=15] (r5) to node[above] {$$} (r12); 	
\draw[midarrow,bend right=15] (r5) to node[above] {$$} (r4); 	
\draw[midarrow,bend right=15] (r5) to node[above] {$$} (r7); 	
\draw[midarrow,bend right=15] (r5) to node[above] {$$} (r8);

\draw[midarrow,bend left=15] (r6) to node[below] {$$} (r9); 	
\draw[midarrow,bend right=15] (r6) to node[above] {$$} (r10); 	
\draw[midarrow,bend right=15] (r6) to node[above] {$$} (r2); 	
\draw[midarrow,bend right=15] (r6) to node[above] {$$} (r3); 	
\draw[midarrow,looseness=15] (r6) to[out=110,in=150] node[left] {$$} (r6); 	
\draw[midarrow,bend right=15] (r6) to node[above] {$$} (r12);

\draw[midarrow,bend right=15] (r7) to node[above] {$$} (r1); 	
\draw[midarrow,bend right=15] (r7) to node[above] {$$} (r5); 	
\draw[midarrow,bend right=15] (r7) to node[above] {$$} (r11); 	
\draw[midarrow,bend right=15] (r7) to node[above] {$$} (r4); 	
\draw[midarrow,bend left=15] (r7) to node[below] {$$} (r8); 	
\draw[midarrow,bend right=15] (r7) to node[above] {$$} (r9);

\draw[midarrow,bend right=15] (r8) to node[above] {$$} (r9); 	
\draw[midarrow,bend right=15] (r8) to node[above] {$$} (r10); 	
\draw[midarrow,bend right=15] (r8) to node[above] {$$} (r2); 	
\draw[midarrow,bend left=15] (r8) to node[below] {$$} (r3); 	
\draw[midarrow,bend right=15] (r8) to node[above] {$$} (r6); 	
\draw[midarrow,bend right=15] (r8) to node[above] {$$} (r12); 	

\draw[midarrow,looseness=15] (r9) to[out=110,in=150] node[left] {$$} (r9); 	
\draw[midarrow,bend left=15] (r9) to node[above] {$$} (r1); 	
\draw[midarrow,bend right=15] (r9) to node[above] {$$} (r2); 	
\draw[midarrow,bend right=15] (r9) to node[above] {$$} (r7); 	
\draw[midarrow,bend right=15] (r9) to node[above] {$$} (r11); 	
\draw[midarrow,bend right=15] (r9) to node[above] {$$} (r12); 	
\draw[midarrow,bend right=15] (r9) to node[above] {$$} (r4); 

\draw[midarrow,bend right=15] (r10) to node[above] {$$} (r5); 
\draw[midarrow,bend right=15] (r10) to node[above] {$$} (r6); 
\draw[midarrow,bend left=15] (r10) to node[above] {$$} (r9); 
\draw[midarrow,looseness=15] (r10) to[out=110,in=150] node[left] {$$} (r10); 
\draw[midarrow,bend right=15] (r10) to node[above] {$$} (r2); 
\draw[midarrow,bend right=15] (r10) to node[above] {$$} (r3); 
\draw[midarrow,bend right=15] (r10) to node[above] {$$} (r6); 
\draw[midarrow,bend right=15] (r10) to node[above] {$$} (r7); 
\draw[midarrow,bend right=15] (r10) to node[above] {$$} (r8); 
\draw[midarrow,bend right=15] (r10) to node[above] {$$} (r12); 

\draw[midarrow,bend right=15] (r11) to node[above] {$$} (r1); 
\draw[midarrow,bend right=15] (r11) to node[above] {$$} (r5); 
\draw[midarrow,looseness=15] (r11) to[out=110,in=150] node[left] {$$} (r11);  
\draw[midarrow,bend right=15] (r11) to node[above] {$$} (r4);  
\draw[midarrow,bend right=15] (r11) to node[above] {$$} (r7);  
\draw[midarrow,bend right=15] (r11) to node[above] {$$} (r8);

\draw[midarrow,bend right=15] (r12) to node[above] {$$} (r9);  
\draw[midarrow,bend right=15] (r12) to node[above] {$$} (r10);  
\draw[midarrow,bend right=15] (r12) to node[above] {$$} (r2);  
\draw[midarrow,bend right=15] (r12) to node[above] {$$} (r3);  
\draw[midarrow,bend right=15] (r12) to node[above] {$$} (r6);  
\draw[midarrow,looseness=15] (r12) to[out=110,in=150] node[left] {$$} (r12);

\end{tikzpicture}
\caption{Sofic graph for $3$-fold rotation\label{fig:soficgraph}}
\end{figure}
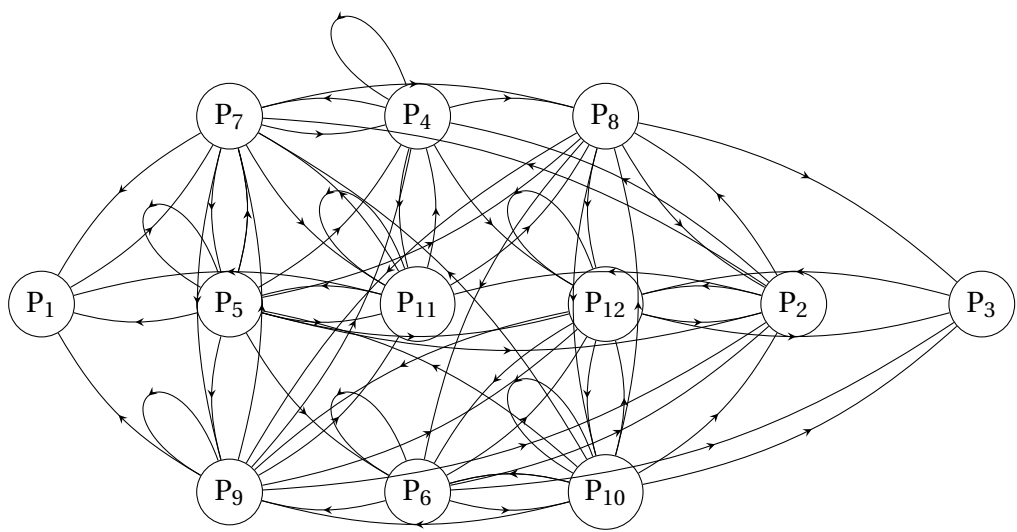
\end{exam}

\begin{exam}
\label{Four}
This example is a kind of a {\it square root system} of
the negative beta expansion introduced by
Ito-Sadahiro 
\cite{Ito-Sadahiro:09}. 
Let $\zeta=i$ and set $\eta_1=1$, $\eta_2=\beta i$ and 
$\xi= -1-\beta i$. We have
$$
T\left(x+yi\right)= 
-\beta y -\left\lfloor -\beta y + 1 \right\rfloor + \beta x i.
$$
By taking its square, we can separate the variables:
$$
T^2\left(x+yi\right)= 
-\beta^2 x - \left\lfloor -\beta^2 x + 1 \right\rfloor + \left(-\beta^2 y -\beta 
\left\lfloor -\beta y + 1 \right\rfloor \right)i.
$$
Thus we can study this map Gaussian coordinate-wise by defining 
$$
f(x)=-\beta^2 x - \left\lfloor -\beta^2 x + 1 \right\rfloor,
$$
a $1$-dimensional piecewise expansive map 
from $\left[-1, 0\right)$ 
to itself and 
$$
g(y)= -\beta^2 y -\beta \left\lfloor -\beta y + 1 
\right\rfloor
$$ 
defined on $\left[-\beta,0\right)$. 
We easily see that $f$ and $g$ give isomorphic systems 
through the relation $g(\beta x) = \beta f(x)$.
Liao-Steiner \cite{Liao-Steiner:12}
showed that the unique ACIM of $f$ is equivalent to 
the $1$-dimensional Lebesgue measure if and only if $\beta^2\ge (1+\sqrt{5})/2$.
Thus the ACIM of $T$ is equivalent to the 
$2$-dimensional Lebesgue measure if and only if $\beta\ge \sqrt{(1+\sqrt{5})/2}$.
In view of the shape of $f$, 
one see that if $\beta^2$ is a Pisot number, then the system 
$(\X, T)$ is sofic (cf. Theorem 3.3 in \cite{Kalle-Steiner:12}).
This give examples of sofic rotational beta expansion 
beyond the scope of Theorem \ref{Sofic}. 
One can also show that when $\beta$ is the 
Salem number whose minimum polynomial is  $x^4-x^3-x^2-x+1$, the system becomes
sofic. 

This example is essentially $1$-dimensional. 
We do not yet succeed in giving a `genuine' $2$-dimensional 
sofic rotational beta expansion beyond Theorem \ref{Sofic}.
\end{exam}

\begin{exam}
\label{Five}
Let $\xi = 0$, $\eta_1=1$ and $\eta_2=\zeta=\exp(2\pi i/5)$. 
Let $\beta = \frac{1+\sqrt{5}}{2}$. 
We describe the symbolic dynamical system associated to given
rotation beta transformation through its sofic graph. 
Here, we use the map $U$ instead of $T$.  
The alphabet $\mathcal{A}=\Delta +\left( \begin{matrix} b_1\\b_2 \\ \end{matrix}\right)=\Delta$ is given by 
$$\left\{
a=\left( \begin{matrix} -2\\0 \\ \end{matrix}\right),
b=\left( \begin{matrix} -1\\0 \\ \end{matrix}\right),
c=\left( \begin{matrix} -2\\1 \\ \end{matrix}\right),
d=\left( \begin{matrix} -1\\1 \\ \end{matrix}\right),
e=\left( \begin{matrix} -2\\2 \\ \end{matrix}\right),
f=\left( \begin{matrix} -1\\2 \\ \end{matrix}\right)
\right\}.$$  

The partition of the fundamental region $[0,1)^2$ is 
given in Figure \ref{fig:fivefold}. 
The sofic graph 
is described in Table \ref{table:five}. 
Since the incidence matrix of this graph
is primitive, we can determine the ACIM whose density is positive and constant on
each partition.
Therefore the ACIM is equivalent to the Lebesgue measure, although
we can not apply Theorem \ref{Ergodic} for $\beta<2$.

\begin{table}[ht]
\caption{} 
\centering 
\begin{tabular}{l l c} 
\hline
\hline 
$n$ &  $U(P_n)$ & $\delta$  \\ [0.5ex]  

\hline 
1 &  $P_{28} \cup  P_{29}$ & $b$ \\ 
2 &  $P_{30} \cup P_{32}\cup P_{33} $ & $b$  \\
3 &  $P_{31}\cup P_{34}\cup P_{35} $  &   $b$ \\
4 &  $P_9 \cup P_{12}\cup P_{19}\cup P_{20}\cup P_{21}\cup P_{22}  $ &   $b$ \\
5 &  $P_6\cup P_{7}\cup P_{18}  $ &   $b$\\ 
6 &  $P_{11}$ & $b$ \\
	&	 $P_2$ & $d$\\
7 &  $P_{37} \cup P_{40}$ & $a$\\
  &  $P_{8} \cup P_{10}$ & $b$\\
  &  $P_{26} \cup P_{27}\cup P_{28}\cup P_{29}$ & $c$\\
	&  $P_{1}$ & $d$\\
8 &  $P_{36}\cup P_{38}\cup P_{39}$ & $a$ \\
  &  $P_{25}$ & $c$ \\
9 &  $P_{30} \cup P_{31}$ & $c$\\
10  &  $P_{23}\cup P_{24}$ & $a$ \\
    &  $P_{13}\cup P_{14}\cup P_{15}\cup P_{16}$ & $c$ \\
11	&  $P_{17}\cup P_{18}$ & $c$ \\
12  &  $P_{19}$ & $c$ \\
13  &  $P_{37}$ & $b$ \\
14  &  $P_{40}$ & $b$ \\
		&  $P_{26}\cup P_{27}$ & $d$ \\
15  &  $P_{36}$ & $b$ \\
16  &  $P_{38}\cup P_{39}$ & $b$ \\
    &  $P_{25}$ & $d$ \\
17  &  $P_{23}\cup P_{24}$ & $b$ \\
		&  $P_{13}\cup P_{14}$ & $d$ \\
18	&  $P_{3}\cup P_{15}\cup P_{16}$ & $d$ \\
19  &  $P_{4}\cup P_{17}$ & $d$ \\
20  &  $P_{33}\cup P_{35}$ & $c$\\
		&  $P_{5}$ & $d$ \\
[1ex]\hline 
\end{tabular}
\qquad
\begin{tabular}{l l c} 
\hline
\hline 
$n$ &  $U(P_n)$& $\delta$  \\ [0.5ex]  
\hline 
21  &  $P_{32}\cup P_{34}$ & $c$ \\
22  &  $P_{20}$ & $c$ \\
23  &  $P_{21}$ & $c$ \\
24  &  $P_{22}$ & $c$ \\
25  &  $P_{28}\cup P_{29}$ & $d$ \\
26  &  $P_{30}\cup P_{32}$ & $d$ \\
27  &  $P_{33}$ & $d$ \\
28  &  $P_{35}$ & $d$ \\
29  &  $P_{31}\cup P_{34}$ & $d$ \\
30  &  $P_{19}\cup P_{20}$ & $d$ \\
31  &  $P_{6}\cup P_{18}$ & $d$ \\
32  &  $P_{9}\cup P_{21}$ & $d$ \\
33  &  $P_{11}\cup P_{12}\cup P_{22}$ & $d$ \\
		&  $P_{2}$ & $f$ \\
34  &  $P_{36}\cup P_{37}$ & $c$ \\
		&  $P_{7}\cup P_{8}$ & $d$ \\
35	&  $P_{39}\cup P_{40}$ & $c$ \\
    &  $P_{10}$ & $d$\\
    &  $P_{27}\cup P_{28}$ & $e$ \\
  	&  $P_{1}$ & $f$\\
36	&  $P_{38}$ & $c$ \\
  	&  $P_{25}\cup P_{26}\cup P_{29}$ & $e$ \\
37	&  $P_{30}\cup P_{31}$ & $e$ \\
38  &  $P_{23}$ & $c$ \\
  	&  $P_{13}\cup P_{15}$ & $e$ \\
39	&  $P_{24}$ & $c$ \\
  	&  $P_{14}\cup P_{16}$ & $e$ \\
40	&  $P_{17}\cup P_{18}\cup P_{19}$ & $e$ \\
&&\\
&&\\
[1ex]\hline 
\end{tabular}
\label{table:five} 
\end{table}
\end{exam}

\begin{figure}[htp]
\centering
\includegraphics[scale=0.11]{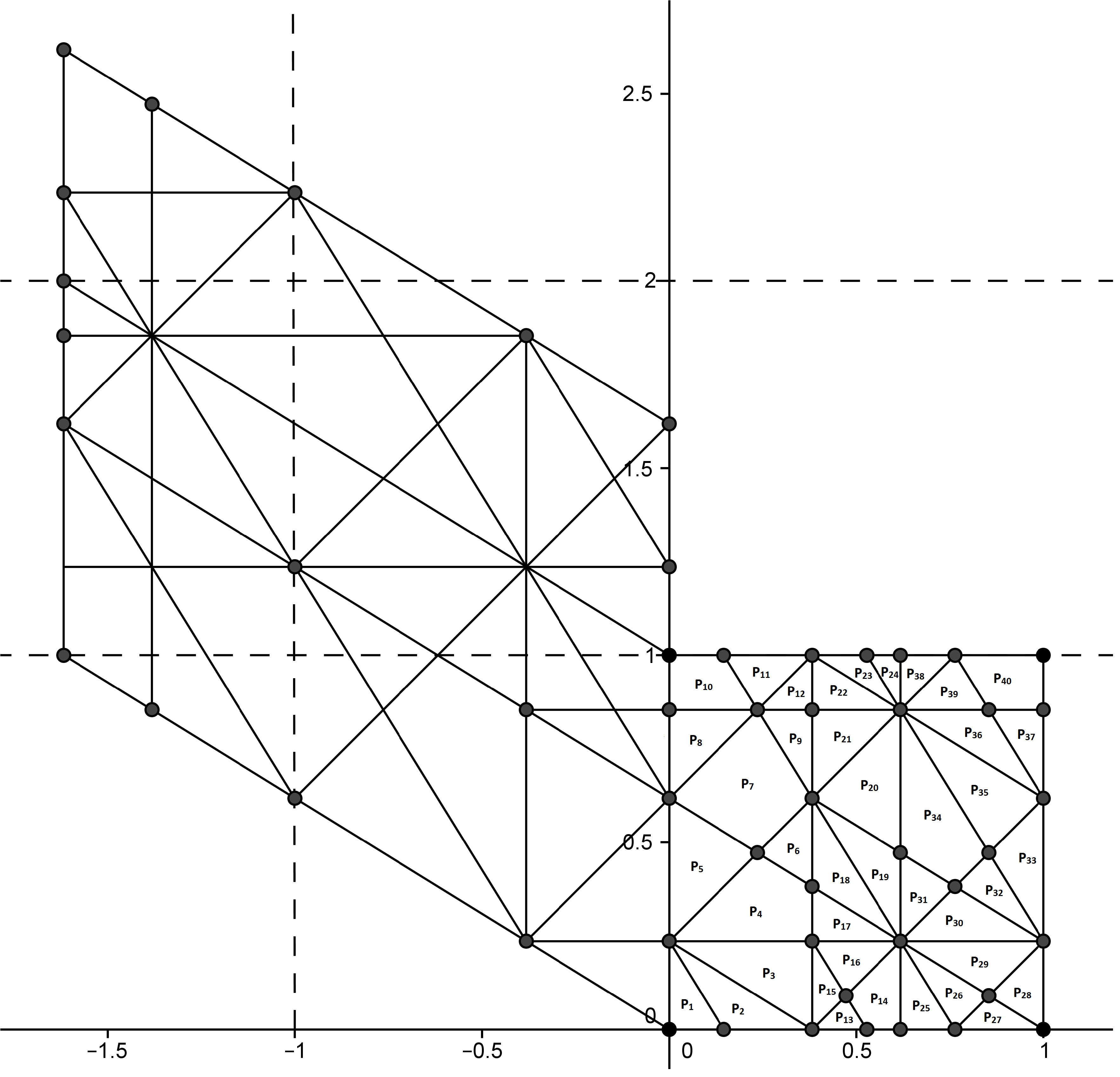}
\caption{5-fold sofic case\label{fig:fivefold}}
\end{figure}

\begin{exam}
\label{Seven}
Let $\xi = 0$, $\eta_1=1$ and $\eta_2=\zeta=\exp(2\pi i/7)$. 
Let $\beta = 1+2\cos (2 \pi /7)\approx 2.24698$, a cubic Pisot number whose minimum 
polynomial is $x^3-2x^2-x+1$.
From $r(\L)=1/(2\cos(\pi/7))$, $w(\X)=\sin(2\pi/7)$ we have $\beta>B_1\approx
2.00272$ and there is a unique ACIM by Theorem \ref{Ergodic}, 
but $\beta<B_2\approx 2.41964$.
From Theorem \ref{Sofic},
we know that the corresponding dynamical system is sofic.
Figure \ref{fig:seven} shows
the sofic dissection of $\X$ by 224 discontinuity segments.
The number of states of the sofic graph is 3292 (!), computed by 
Euler's formula.
It is possible to show that
the corresponding incidence matrix of the sofic 
graph is primitive, and consequently
the ACIM is equivalent to the Lebesgue measure.

\begin{figure}[htp]{}
	\centering
		\includegraphics[width=\textwidth]{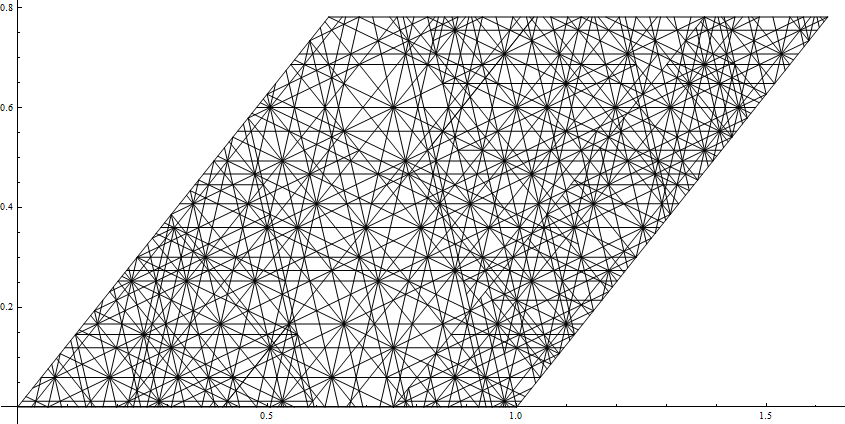}
\caption{Sofic 7-fold rotation\label{fig:seven}}
	\end{figure}
\end{exam}

\end{section}


\end{document}